# LOCAL LIMIT THEORY AND LARGE DEVIATIONS FOR SUPERCRITICAL BRANCHING PROCESSES

BY PETER E. NEY AND ANAND N. VIDYASHANKAR

*University of Wisconsin and University of Georgia*

In this paper we study several aspects of the growth of a supercritical Galton–Watson process $\{Z_n : n \geq 1\}$, and bring out some criticality phenomena determined by the Schröder constant. We develop the local limit theory of $Z_n$, that is, the behavior of $P(Z_n = v_n)$ as $v_n \nearrow \infty$, and use this to study conditional large deviations of $\{Y_{Z_n} : n \geq 1\}$, where $Y_n$ satisfies an LDP, particularly of $\{Z_n^{-1} Z_{n+1} : n \geq 1\}$ conditioned on $Z_n \geq v_n$.

**1. Introduction.** In this paper we study the large deviations of a "random average" indexed by a supercritical branching process and related aspects of the growth rate of the branching process. We introduce "conditional large deviation theory" and establish that certain functionals based on the $n$th generation population size satisfy the conditional large deviation principle. In the process we also establish a phase transition in the rate of growth of the branching process based on the values of a parameter $\alpha$, the so-called Schröder constant. The main technical tool is a local limit theorem which substantially unifies, sharpens and extends the existing results in the literature.

We begin by considering the single type Galton–Watson branching process $\{Z_n, n \geq 1\}$ initiated by a single ancestor, that is, $Z_0 \equiv 1$. We denote the offspring distribution by $\{p_j : j \geq 0\}$, the mean of the offspring distribution by $m$ $(> 1)$, and the probability generating function by $f(s)$ for $0 \leq s \leq 1$, that is,

$$P[Z_1 = j] = p_j, \qquad j \geq 0, \qquad f(s) = \sum_{j \geq 1} s^j p_j \quad \text{and} \quad m \equiv \sum_{j \geq 1} j p_j.$$

Let $\gamma = f'(q)$, where $q = P(Z_n = 0 \text{ for some } n \geq 1)$ is the extinction probability. Let $\{\xi_{n,i} : i \geq 1, n \geq 1\}$ be i.i.d. random variables with $P(\xi_{n,i} = k) = p_k$, and interpret them as the number of offspring of the $i$th parent in the $n$th generation. Then $Z_{n+1} = \sum_{k=1}^{Z_n} \xi_{n,k}$.

The focus of the paper is on the the large deviation behavior of the ratio $R_n = Z_n^{-1} Z_{n+1}$ and some of its generalizations. By the branching property, this ratio can be expressed as

$$R_n = \frac{1}{Z_n} \sum_{i=1}^{Z_n} \xi_{n,i} = \frac{1}{Z_n} \sum_{i=1}^{Z_n} \xi_i \equiv \overline{\xi}_{Z_n},$$

where $\{\xi_i, i \geq 1\}$ are i.i.d. with $P(\xi_1 = k) = p_k$, $k \geq 0$. Its behavior sheds light on the evolution of branching populations and is also important in statistical inference for branching processes since it is the maximum likelihood estimator of the mean $m$ when $(Z_n, Z_{n+1})$ are observed. It is thus relevant to study its rate of convergence; doing so, via Bahadur efficiency, leads to questions about large deviations of $R_n$. Furthermore, $R_n$ being a random average of i.i.d. random variables, it is natural to enquire whether it has large deviation properties along the lines of Cramér's theorem.

Large deviations of $\{R_n : n \geq 1\}$ were previously investigated by Athreya [1] and Athreya and Vidyashankar [4]. It was established in those papers that, under an exponential moment hypothesis, if $a > m$ and $\gamma > 0$, then

(1) $$0 < \lim_{n \to \infty} \frac{1}{n} \log P(R_n > a | Z_n > 0) = -\log \gamma.$$

Thus the sequence $R_n$ satisfies the LDP (see [10] for definition) with the constant rate function $I(x) = -\log \gamma$. The limit is "degenerate" in the sense that the rate function is independent of $a$.

It is also puzzling, at first sight, that though $R_n = Z_n^{-1} \sum_{i=1}^{Z_n} \xi_i$ is a "sample mean," the Cramér rate function $I = \Lambda^\star$ (the convex conjugate of the logarithmic generating function of $\xi_1$) does not appear in the rate. A heuristic explanation of this behavior comes from the fact that the contributions to $P(R_n \geq a)$ come from small values of $Z_n$. More precisely, a straightforward argument shows that

(2) $$\lim_{n \to \infty} P(Z_n = k | R_n \geq a) = a(k),$$

where $\{a(k) : k \geq 1\}$ is a probability distribution. Thus, even though $Z_n$ grows like $m^n$ when conditioned on $R_n > a$, $Z_n$ stays small in the sense that it converges (conditionally, in probability) to a proper random variable; that is, the "large" deviations of $R_n$ favor "small" values of $Z_n$. Thus if one forces $Z_n$ to be large, then it is conceivable that $R_n$ may exhibit a "typical" large deviation behavior as described by Cramér's theory. This motivates studying the large deviations of $\{R_n : n \geq 1\}$ by conditioning on $Z_n \geq v_n$ where $v_n \nearrow \infty$. A simple Bayes formula for the conditional probability then leads to studying in detail the behavior of $P(Z_n = v_n)$ for various choices of $v_n$.

Conditioning as a technique to improve the precision of confidence intervals has long been used in the statistical literature. It has been argued by Efron, Hinkley and others (see [16, 26, 27]) that when one seeks to construct confidence intervals for an unknown parameter it should be conditioned on the observed Fisher information. Motivated by this, Sweeting [28] considered the problem of maximum likelihood estimation of the offspring mean when the underlying offspring distribution is geometric. Indeed, Sweeting establishes that the statistic $T_n \equiv (\sum_{j=0}^{n} Z_{j-1})^{-1} \sum_{j=0}^{n} Z_j$, conditioned on $\sum_{j=0}^{n} Z_{j-1} = V_n$, where $V_n \sim m^n$, appropriately centered and normalized converges to a Gaussian random variable with mean 0 and variance equal to the Fisher information. Furthermore, Sweeting [29] demonstrates the improved performance of these conditional confidence intervals through simulations. For further work on conditioning in the context of branching processes, see [6, 18, 17]. We of course focus on the statistic $R_n$. The observed conditional Fisher information in this case is $\sigma^{-2} Z_n$.

A further motivation for the present study comes from a class of conditioned limit laws similar to Gibbs conditioning. In its simplest form the latter describes the behavior of the conditional distributions

$$P\bigg\{(X_1,\ldots,X_k) \in \cdot \Big| \frac{1}{n}\sum_{i=1}^{n} X_i \in A\bigg\},$$

where $X_1, X_2, \ldots$ are i.i.d. random variables and $A$ is a Borel set. Extension of this concept to branching processes suggests studying the branching distributions

(3) $$P\bigg\{(Z_1,\ldots,Z_k) \in \cdot \Big| \frac{1}{Z_n}\sum_{i=1}^{Z_n} \xi_{n,i} \in A\bigg\}.$$

Since the conditioning event is $\{R_n \in A\}$, or more generally can be taken of the form $\{(Z_n, Z_{n+1}) \in A \times B\}$, behavior of (3) yields information about the *past* structure of the branching population, based on the information about its *present* (namely $Z_n$ and $Z_{n+1}$). This formulation has practical significance, for example, in the area of molecular evolution [9, 10, 23, 30]. Now the analysis of (3) requires a careful study of the large deviations of $\{R_n : n \geq 1\}$ which is carried out in the present paper. The behavior of (3) itself requires additional technical tools which are developed and treated separately in another paper.

The rest of the paper is organized as follows: Section 2 contains a summary of results and related discussions, while Section 3 contains proofs. Section 4 is devoted to some concluding remarks.

**2. Summary of results and related discussions.** In this section we state the main results of the paper. We begin by describing the relevant branching process background so as to ease the discussion and exhibit the significance of our results.

2.1. *Branching process background.* Let $\{Z_n : n \geq 1\}$ be a single type supercritical branching process with $Z_0 \equiv 1$ and mean $m$. The sequence $\{W_n \equiv m^{-n} Z_n ; n \geq 1\}$ plays an important role in the study of the limit theory of supercritical branching processes. It is well known that $\{W_n : n \geq 1\}$ is a nonnegative martingale sequence and, hence, converges with probability one to a nonnegative random variable $W$; under the further assumption that $E(Z_1 \log Z_1) < \infty$, the limit random variable is nontrivial with an absolutely continuous density $w(\cdot)$, except for a possible atom at 0 (see [3] for details). The behavior of the density near 0 has been investigated by Dubuc [13] who showed that for $0 < x < 1$ there exists universal constants $0 < C_1 < C_2 < \infty$ such that

$$(4) \qquad C_1 x^{\alpha-1} \leq w(x) \leq C_2 x^{\alpha-1},$$

where $\alpha$ is a solution to $\gamma = m^{-\alpha}$ and is assumed to be finite. The finiteness of $\alpha$ is equivalent to $p_0 + p_1 > 0$ and, drawing on the functional iterations literature, offspring distributions with this property are said to belong to the *Schröder case*. A further refinement of the above estimate has been considered by Biggins and Bingham [7]. The process with $\alpha = \infty$ grows exponentially fast at all times and as such has a different probabilistic structure and is frequently referred to as the *Böttcher case*.

The quantity $\alpha$ shows up in several deep results in the theory of supercritical branching processes and will play a critical role in our study as well. Karlin and McGregor [21, 22] studied the problem of embeddability of discrete-time branching processes into continuous-time branching processes. The Karlin–McGregor function

$$K(s) \equiv s^\alpha Q(\phi(s)) \qquad \text{where } \phi(s) \equiv E(e^{-sW})$$

and its constancy has been the subject of much study (see [14] and the references therein). For the definition of $Q(\cdot)$ see below. Karlin and McGregor conjectured that $K(\cdot)$ is constant exactly when the the discrete-time process is embeddable into a continuous-time process. Building on their work, Dubuc [13] established that a discrete-time branching process [satisfying $E(Z_1 \log Z_1) < \infty$] is embeddable into continuous-time Markov branching process if and only if

$$\lim_{x \searrow 0} w(x) x^{1-\alpha} \qquad \text{exists and is finite.}$$

He further established that this is equivalent to the existence of the limit as $s \nearrow 1$ of $Q(s)(1-s)^\alpha$, where for $0 \leq s \leq 1$,

$$(5) \qquad Q(s) \equiv \lim_{n \to \infty} Q_n(s) \equiv \lim_{n \to \infty} \frac{f_n(s) - q}{\gamma^n}.$$

The function $Q(\cdot)$ satisfies the functional equation $Q(f(s)) = \gamma Q(s)$ and $Q(q) = 0$ and $Q(1) = \infty$; it can be extended to an analytic function in the interior of the unit disc and hence has the power series representation

$$Q(s) = \sum_{j \geq 0} q_j s^j \qquad \text{for } 0 \leq s < 1. \tag{6}$$

In the work on large deviations of branching processes (and in other contexts as well; see [19, 20]) it has been established [1, 24] that the integrability of $Q(\cdot)$ near 1, that is, finiteness of $I \equiv \int_0^1 Q(s)\,ds$, plays an important role. In the embeddable situation $0 < \alpha < 1$ and the above-mentioned result of Dubuc readily yields that $I < \infty$. Of course, one can establish the finiteness of $I$, under the assumption $\alpha < 1$, without invoking any embeddability issues [1, 24, 25].

In this paper we make a detailed study of $P(Z_n = v_n)$ for various ranges of values of $v_n$. This problem was studied in the 1970s for the extreme cases when $v_n \sim m^n$ for a complicated restriction on the range of values of $\alpha$ or when $v_n$ is a constant [2, 13, 14]. The behavior of $P(Z_n = v_n)$ for the entire range of values of $v_n = \mathbf{O}(m^n)$ and the entire range of values of $\alpha$ has been open since that time. We will present a unified solution to this problem and show that this rate is intimately connected to the rate of convergence of $I_n = \int_c^1 Q_n(s)\,ds$ to $I = \int_c^1 Q(s)\,ds$, where $0 < c < 1$.

2.2. *Local limit theorem.* We begin this section by recalling the known local limit results. For the sake of simplicity, we will assume that $p_0 = 0$ (see Remark 4). When $v_n$ is a constant, one can use the analyticity of $Q(\cdot)$ in the interior of the unit disc to show that

$$\lim_{n \to \infty} \frac{P(Z_n = j)}{p_1^n} = q_j, \tag{7}$$

where $q_j$'s are defined in (3). A natural question is the behavior of $P(Z_n = v_n)$ when $v_n$ grows with $n$. If one assumes $E(Z_1^2) < \infty$, then Athreya and Ney [2] (see also [13]) have established the rate of convergence to 0 of the difference $|m^n P(Z_n = j) - w(jm^{-n})|$, when $j(n) \sim m^n$. Assuming only the finiteness of the mean, Dubuc and Seneta [15] proved a weaker form of the above result. It is unclear from these works if $\alpha$ plays any role in the behavior of $P(Z_n = v_n)$. For the early history of the local limit theorem, see the paper of Athreya and Ney [2].

Our first result is a new local limit theorem which provides the asymptotic behavior of $P(Z_n = v_n)$ in the entire range $v_n = \mathbf{O}(m^m)$, and demonstrates the critical role of $\alpha$. The theorem covers both the Schröder and Böttcher cases.

THEOREM 1. *Assume that $E(Z_1 \log Z_1) < \infty$. Let $\{v_n : n \geq 1\}$ be a sequence of integers such that $v_n \nearrow \infty$ as $n \to \infty$ and $v_n = \mathbf{O}(m^n)$. Then there exists constants $0 < C_1 < C_2 < \infty$ such that*

$$C_1 \leq \liminf_{n \to \infty} \frac{P(Z_n = v_n)}{A_n} \leq \limsup_{n \to \infty} \frac{P(Z_n = v_n)}{A_n} \leq C_2,$$

*where*

$$A_n = \begin{cases} p_1^n v_n^{\alpha - 1}, & \text{if } \alpha < 1, \\ k_n p_1^n, & \text{if } \alpha = 1, \\ m^{-n}, & \text{if } \alpha > 1 \text{ (possibly } \infty\text{)}, \end{cases}$$

*and $k_n = [n - \frac{\log v_n}{\log m} + 1]$. Furthermore, if $v_n = m^{n-k_n}$ for some sequence of integers $k_n = \mathbf{O}(n)$ as $n \to \infty$, then $\lim_{n \to \infty} A_n^{-1} P(Z_n = v_n)$ exists and is positive and finite.*

REMARK 1. When $E(Z_1 \log Z_1) = \infty$ but $m < \infty$, the rate of convergence will depend on the Senata constants. The extension of Theorem 1 to this general case is difficult and will be taken up separately.

REMARK 2. Observe that when $n$ is fixed and $\alpha < 1$, $P(Z_n = v_n)$ increases if $v_n$ decreases; thus, roughly speaking, $P(Z_n = z)$ favors smaller values of $z$. This phenomenon will reappear in different guises later. It is further suggested by the fact that from (1) [when $E(Z_1 \log Z_1) < \infty$], the mode of the distribution of $W$ is closer to 0. The case $\alpha = 1$ is somewhat surprising, since such an argument using the density does not hold.

The following uniform estimate on the behavior of $\sup_{j \geq v_n} P(Z_n = j)$ is of interest in its own right and will be used in the proof of our large deviation theorems. The result covers both the Schröder and Böttcher cases.

THEOREM 2. *Under the conditions of Theorem 1, the following prevail:*

1. *If $\alpha < 1$, then there exist universal constants $0 < C_1 < C_2 < \infty$ such that*

(8) $$C_1 p_1^n \leq \sup_{j \geq v_n} j^{\alpha - 1} P(Z_n = j) \leq C_2 p_1^n.$$

2. *If $\alpha = 1$, then there exist universal constants $0 < C_3 < C_4 < \infty$ such that*

(9) $$C_3 p_1^n k_n \leq \sup_{j \geq v_n} P(Z_n = j) \leq C_4 p_1^n k_n,$$

*where $k_n = [n - \frac{\log v_n}{\log m}]$ ($[x]$ refers to the largest integer in $x$).*

3. *If $\alpha > 1$ (possibly $\infty$), then there exist universal constants $0 < C_5 < C_6 < \infty$ such that*

(10) $$C_5 m^{-n} \leq \sup_{j \geq v_n} P(Z_n = j) \leq C_6 m^{-n}.$$

REMARK 3. Theorems 1 and 2 in the Böttcher case, that is, $j_0 \equiv \inf\{j : p_j > 0\} \geq 2$, are interesting since the process has a different probabilistic behavior. Of course, $P(Z_n = v_n) = 0$, whenever $v_n < j_0^n$. Thus the local limit theorem yields nontrivial results when $j_0^n < v_n$ and $v_n = \mathbf{O}(m^n)$. Dubuc [14] has given a detailed analysis of the moment generating function of $W$ and Bingham [8], using a Tauberian argument, has elucidated the behavior of $P(W \leq x)$ as $x \downarrow 0$. A crucial role is played by the parameter $\beta[\equiv (\log m)^{-1} \log j_0]$ that relates the minimum family size and the mean of the offspring distribution. By definition, $0 < \beta < 1$. Then the result of Bingham [8] shows that as $x \downarrow 0$,

$$(11) \qquad -\log P(W \leq x) \sim \frac{\tau}{x^{\beta/(1-\beta)}},$$

where $\tau$ is an interesting constant. The above result can be used to show that the density $w(\cdot)$ of the random varible $W$ decays to 0 as $x$ decreases to 0 exponentially fast in a manner dictated by a function of $\beta$. Thus it is reasonable to expect that the behavior of $P(Z_n = v_n)$, just as in the Schröder case, depends on the range of values of $\beta$. It turns out, however, that since the Böttcher case corresponds to a situation when $\alpha = \infty$, the rate is just $m^{-n}$ (Case $\alpha > 1$) and independent of $\beta$.

REMARK 4. If $p_0 > 0$, then Theorems 1 and 2 will hold with $p_1$ replaced by $\gamma = f'(q)$.

2.3. *Large deviations.* In this section we state our large deviation results. Since $R_n$ is the mean of a random sum, one is motivated to ask more generally about the large deviations behavior of $S(N_n)$ where $S(n) = \sum_{1=1}^{n} X_i$ where $\{X_i : i \geq 1\}$ is a sequence of i.i.d. random variables and $N_n$ ($\nearrow \infty$ in probability) is a sequence of random variables. When one normalizes the sequence $\{S(N_n) : n \geq 1\}$ by $n$ then the large deviation behavior of $n^{-1}S(N_n)$ is easily established [see [11], Exercise 2.3.19(b)]. In fact, the conclusion is that $\{n^{-1}S(N_n) : n \geq 1\}$ satisfies an LDP with the rate function being the Legendre–Fenchel transform of $\Lambda(\Lambda_X(\theta))$, where $\Lambda_X(\theta)$ is the logarithmic moment generating function of $X$, and for some $n \to \infty$,

$$(12) \qquad \Lambda(\lambda) = \lim_{n \to \infty} \frac{1}{n} \log E(e^{\lambda N_n}),$$

where the above limit is assumed to exist in a neighborhood of the origin. The situation is quite different, however, when one normalizes the sequence $\{S(N_n) : n \geq 1\}$ by $N_n$ itself, which is the case with $R_n$. Our first proposition considers this problem in a slightly more general context.

PROPOSITION 1. *Let $\{Y(n) : n \geq 1\}$ be a sequence of random variables that satisfy the LDP with the rate function $I_Y(\cdot)$ and speed (normalizing* sequence) $r(n)$ and $\{N_n : n \geq 1\}$ denote a sequence of integer valued random variables such that $N_n \nearrow \infty$; assume further that for $a_n \to \infty$ and $\theta > 0$,

$$\lim_{n \to \infty} \frac{1}{a_n} \log E(e^{-r(N_n)\theta}) = \Lambda(\theta). \tag{13}$$

Let

$$\mu_n(\cdot) = P(Y(N_n) \in \cdot). \tag{14}$$

Then the sequence of measures $\{\mu_n : n \geq 1\}$ satisfy the LDP with rate function $-\Lambda(I_Y(\cdot))$ and speed $a_n$.

The point of the proposition is that the rate function depends explicitly on the behavior of $\Lambda(\cdot)$. Indeed, if $\Lambda$ is a constant, then the rate function is degenerate in the sense that the rate is "independent" of the set for which the large deviation is studied. We have seen that this general phenomenon is exhibited when one considers the large deviations of $R_n$.

Thus to bring out the dependence of the rate on the set being considered, we condition on $N_n \geq v_n$ or more generally, $N_{n-k} \geq v_{n-k}$. The following form of Proposition 1 brings out the ingredients necessary for the conditional analysis in Theorem 3.

PROPOSITION 2.  *Let $\{Y(n) : n \geq 1\}$ denote a sequence of random variables that satisfy the LDP with rate function $I_Y(\cdot)$ and speed $r(n)$, and $\{N_n : n \geq 1\}$ denote a sequence of integer valued random variables such that $N_n \nearrow \infty$ and is independent of $Y_n$. Let*

$$\mu_n(\cdot) = P(Y(N_n) \in \cdot | N_{n-k} \geq v_{n-k}). \tag{15}$$

*Let $\{N_n^1(k) : n \geq 1\}$ denote the sequence of random variables with distribution*

$$P(N_n^1(k) = j) \equiv P(N_n = j | N_{n-k} \geq v_{n-k}). \tag{16}$$

*Assume that there exists $a_n \to \infty$ such that, for all $\theta > 0$,*

$$\lim_{n \to \infty} \frac{1}{a_n} \log E(\exp(-\theta r(N_n^1(k)))) \equiv K(\theta) \tag{17}$$

*and the limit is continuous. Then the sequence of measures $\{\mu_n : n \geq 1\}$ satisfies the LDP with rate function $-K(I_Y(\cdot))$ and speed $a_n$.*

Of course, in order for the above proposition to have real substance in particular cases, the speed sequence $\{a_n\}$ and the limit $K(\theta)$ in (16) must be explicitly determined. Indeed when $N_n = Z_n$, we can make these determinations using the local limit theory and the results are stated in the following theorem.

THEOREM 3. *Let $\{Y(i) : i \geq 1\}$ be any sequence of random variables that satisfies the LDP with a "good" rate function $I(\cdot)$ and speed $n$. Assume that $E(Z_1 \log Z_1) < \infty$ and $0 < \alpha < \infty$. Let*

(18) $\qquad \mu_{n,k}(A) = P(Y(Z_n) \in A | Z_{n-k} \geq v_{n-k}), \qquad n \geq k, A \subset R,$

*where $\{v_n : n \geq 1\}$ is a sequence of positive integers increasing to infinity.*

1. *If $\lim_{n \to \infty} n^{-1} v_{n-k} = b, 0 \leq b < \infty$, then the sequence $\mu_{n,k}$, where $k$ is fixed, satisfies the LDP with rate $\tilde{I}(x)$ and speed $v_{n-k}$, where*

(19) $\qquad\qquad\qquad \tilde{I}(x) = -\log f_k(e^{-I(x)}) + bB$

   *and $B = -\log p_1$ if $\alpha \leq 1$ while $B = \log m$ if $1 < \alpha < \infty$.*
2. *If $b = 0$, $\tilde{I}$ reduces to $-\log f_k(e^{-I(x)})$.*
3. *If $\lim_{n \to \infty} v_{n-k}^{-1} n \to \infty$, then the sequence $\mu_{n,k}$, where $k$ is fixed, satisfies the LDP with constant rate function $B$ and speed $n$. Furthermore, $\tilde{I}(\cdot)$ is a good rate function.*

REMARK 5. The trichotomy $b = 0$, $0 < b < \infty$, $b = \infty$ shows the critical role of the rate $v_n = n$; namely whether $v_n$ grows faster or slower than $n$. The limit described above may not exist if $\lim_{n \to \infty} v_{n-k}^{-1} n$ does not exist. In fact, one can construct a sequence $v_n$ such that the limit is different along different subsequences. However, in these situations the large deviations upper bound and lower bound will still hold.

The above theorem when specialized to the branching case yields the large deviations for $R_n$. We have the following corollary.

COROLLARY 1. *Assume $E(e^{\theta Z_1}) < \infty$, $\lim_{n \to \infty} n^{-1} v_{n-k} = b$ and $0 < \alpha < \infty$. Then sequence of measures*

(20) $\qquad\qquad \mu_{n,k}^{(1)} \equiv P(R_n \in A | Z_{n-k} \geq v_{n-k}), \qquad n > k$

*satisfies the LDP with the good rate function $\tilde{I}_1(\cdot)$ and speed $v_{n-k}$, where*

(21) $\qquad\qquad\qquad \tilde{I}_1(x) = -\log f_k(e^{-\Lambda^\star(x)}) + bB;$

*here $\Lambda^\star(\cdot)$ is the convex conjugate of the cumulant generating function of the random variable $Z_1$ and $B$ is as in the theorem. If $n^{-1} v_{n-k} \to \infty$, then $\tilde{I}_1(x) = B$.*

REMARK 6. Under an exponential moment condition on $Z_1$, a sharper result along the lines of [5] can be obtained.

Theorem 3 also sheds light on the large deviations of $W_n^{-1} W$. Namely, we have the following corollary.

COROLLARY 2. *Assume* $E(e^{\theta Z_1}) < \infty$, $\lim_{n\to\infty} n^{-1} v_{n-k} = b$ *and* $0 < \alpha < \infty$. *Then the sequence of measures*

$$(22) \qquad \mu_{n,k}^{(2)} \equiv P\left(\frac{W}{W_n} \in A \Big| Z_{n-k} \geq v_{n-k}\right), \qquad n > k$$

*satisfies the LDP with the good rate function* $\tilde{I}_2(\cdot)$ *and speed* $v_{n-k}$, *where*

$$(23) \qquad \tilde{I}_2(x) = -\log f_k(e^{-\Lambda_W^\star(x)}) + bB,$$

*where* $\Lambda_W^\star(\cdot)$ *is the Legendre–Fenchel transform of the cumulant generating function of the random variable* $W$, *and* $B$ *is as in the theorem. If* $n^{-1} v_{n-k} \to \infty$ *then* $\tilde{I}_2(x) = B$.

An upper bound for the rate of convergence of $|W - W_n|$ was obtained by Athreya [1].

Our next theorem treats $\mu_{n,k}$ in (18) when $\alpha = \infty$ (the Böttcher case).

THEOREM 4. *Let* $\{X(i) : i \geq 1\}$ *be any sequence of random variables that satisfy the LDP with a "good" rate function* $I(\cdot)$. *Assume that* $E(Z_1 \log Z_1) < \infty$. *Let*

$$(24) \qquad \mu_{n,k}(A) = P(X(Z_n) \in A | Z_{n-k} \geq v_{n-k}), \qquad n \geq k, A \subset R,$$

*where* $\{v_n : n \geq 1\}$ *is a sequence of positive integers increasing to infinity. Let* $j_0 = \inf\{j \geq 2 : p_j > 0\}$. *If* $\lim_{n\to\infty} v_n^{-1} j_0^n = b$ $(> 0)$ *then the sequence* $\mu_{n,k}$ *(where* $k$ *is fixed) satisfies the LDP with the good rate function* $\tilde{I}_3(x)$ *and speed* $v_{n-k}$, *where*

$$(25) \qquad \tilde{I}_3(x) = -bG(f_k(e^{-I(x)}));$$

*here, for* $0 \leq s \leq 1$,

$$(26) \qquad G(s) = \log s + \log p_{j_0} + \sum_{j \geq 0} \frac{1}{j_0^{j+1}} \log\left(1 + \frac{1 - p_{j_0}}{p_{j_0}} g(f_j(s))\right)$$

*and*

$$(27) \qquad g(s) = \frac{1}{1 - p_{j_0}} \sum_{j \geq j_0 + 1} p_j s^{j - j_0}.$$

*If* $b = 0$ *and* $v_n = \mathbf{O}(m^n)$ *(again for fixed* $k$*) then the sequence* $\mu_{n,k}$, *where* $k$ *is fixed, satisfies the LDP with the good rate function* $\tilde{I}_3(x)$ *and speed* $v_{n-k}$, *where*

$$(28) \qquad \tilde{I}_3(x) = -\log f_k(e^{-I(x)}).$$

A process valued result along the lines of Mogulski's theorem [11] can also be obtained in our conditional setting. We state this as Corollary 3 for case $\alpha < \infty$. Let

$$(29) \qquad R_n(t) \equiv \frac{1}{Z_n} \sum_{i=1}^{[tZ_n]} \xi_{n,i}, \qquad 0 \le t \le 1.$$

When one views $R_n(t)$ as a random variable taking values in $L_\infty[0,1]$, one has the LDP for the sequence of measures $\mu_{n,k}(\cdot) \equiv P(R_n(t) \in \cdot | Z_{n-k} \ge v_{n-k})$ in $L_\infty(0,\infty)$.

COROLLARY 3. *The sequence of measures $\mu_{n,k}$ satisfies the LDP on $L_\infty[0,1]$ with a good rate function*

$$\bar{I}(\phi) = \begin{cases} \int_0^1 I_1(\phi'(t))\, dt, & \text{if } \phi \in \mathcal{AC}, \phi(0) = 0, \\ \infty, & \text{otherwise,} \end{cases}$$

*where $\mathcal{AC}$ denotes the space of absolutely continuous functions and*

$$(30) \qquad I_3(x) = -\log f_k(e^{-\Lambda^*(x)}) - bB$$

*and $B$ is as in the Theorem 3 and $\Lambda^\star(\cdot)$ is as in Corollary 1.*

REMARK 7. A result similar to Corollary 3 also holds for the Böttcher case.

**3. Proofs.** In this section we provide the proofs of results given in Section 2. We begin by recalling that $f_k(s) \equiv E(s^{Z_k})$ is the $k$th iterate of $f(s)$. Let $\psi_k(s) \equiv E(e^{isW_k}) = f_k(e^{ism^{-k}})$ and $\psi(u) = E(e^{iuW})$ denote the characteristic functions of $W_k$ and $W$, respectively. By the inversion theorem (see [3], page 81),

$$(31) \qquad P(Z_n = v_n) = \frac{1}{2\pi} \int_{-\pi}^{\pi} e^{-i\theta v_n} f_n(e^{i\theta})\, d\theta.$$

Making a change of variable $\theta v_n \to u$ we have that

$$(32) \qquad P(Z_n = v_n) = \frac{1}{2\pi v_n} \int_{-\pi v_n}^{\pi v_n} e^{-iu} f_n(e^{iu/v_n})\, du.$$

Since $v_n \le m^n$, there exists $\tilde{k}_n > 0$ such that $v_n = m^{n-\tilde{k}_n}$. Let $k_n = [\tilde{k}_n]$, $a_n = n - k_n$ and $\eta(n) = m^{\tilde{k}_n - k_n}$. Then using functional iteration, we have that

$$(33) \qquad P(Z_n = v_n) = \frac{1}{2\pi v_n} \int_{-\pi v_n}^{\pi v_n} e^{-iu} f_{k_n}(\psi_{a_n}(u\eta(n)))\, du$$

$$= \frac{1}{2\pi v_n} \int_{-\pi v_n}^{0} e^{-iu} f_{k_n}(\psi_{a_n}(u\eta(n))) \, du$$

$$+ \frac{1}{2\pi v_n} \int_{0}^{\pi v_n} e^{-iu} f_{k_n}(\psi_{a_n}(u\eta(n))) \, du$$

(34) $$\equiv J_n(1) + J_n(2).$$

We make a detailed analysis of $J_n(2)$. Note that

(35) $$J_n(2) = \frac{1}{2\pi v_n} \int_{0}^{\pi} e^{-iu} f_{k_n}(\psi_{a_n}(u\eta(n))) \, du$$

$$+ \frac{1}{2\pi v_n} \int_{\pi}^{\pi v_n} e^{-iu} f_{k_n}(\psi_{a_n}(u\eta(n))) \, du$$

(36) $$\equiv J_n(2,1) + J_n(2,2).$$

We first state a "decomposition lemma" that will help us deal with $J_n(2,1)$.

LEMMA 1.   *For any positive integers $r, s$ ($s < r$) the following decomposition holds:*

(37) $$\int_{0+}^{\pi} e^{-iu} f_r(\psi_s(u)) \, du = m^{-r} I^{(1)}(r,s) + \sum_{l=0}^{r} m^{-l} I_l^{(2)}(r,s),$$

*where*

(38) $$I^{(1)}(r,s) \equiv \int_{0}^{\pi m^{-1}} e^{-ium^{-r}} \psi_{s+r}(u) \, du$$

*and*

(39) $$I_l^{(2)}(r,s) \equiv \int_{\pi m^{-1}}^{\pi} e^{-ium^{-l}} f_{r-l}(\psi_{s+l}(u)) \, du.$$

PROOF.   We first decompose the integral on the left-hand side into two components, one over the interval $(0+, \pi m^{-(r+1)})$ and then over the interval $(\pi m^{-(r+1)}, \pi)$. Let

(40) $$J_{r,s}^1 = \int_{0+}^{\pi m^{-(r+1)}} e^{-iu} f_r(\psi_s(u)) \, du$$

and

(41) $$J_{r,s}^2 = \int_{\pi m^{-(r+1)}}^{\pi} e^{-iu} f_r(\psi_s(u)) \, du.$$

Now making a change of variable $u \longmapsto um^{-r}$ and using the definition of $\psi_{r+s}(u)$, we get that

(42) $$J_{r,s}^1 = m^{-r} I^{(1)}(r,s).$$

We can further decompose $J_{r,s}^2$ as

$$(43) \qquad J_{r,s}^2 = \sum_{l=0}^{r} \int_{\pi m^{-(l+1)}}^{\pi m^{-(l)}} e^{-iu} f_r(\psi_s(u))\, du.$$

Now using the fact that

$$(44) \qquad f_r(\psi_s(u)) = f_{r-l}(\psi_{s+l}(um^l))$$

and the change of variable $um^l \longmapsto u$, it follows that

$$(45) \qquad \int_{\pi m^{-(l+1)}}^{\pi m^{-l}} e^{-iu} f_r(\psi_s(u))\, du = m^{-l} I_l^{(2)}(r,s). \qquad \square$$

The next lemma provides a uniform estimate on the characteristic function of the random variable $W_r$. The proof is similar to the proof of Lemma 2 in [15] and hence is omitted.

LEMMA 2. *Let $\{\eta(k): k \geq 1\}$ be a bounded sequence of positive numbers such that $\inf_{n \geq 1} \eta(n) > 0$. Then for every $\varepsilon > 0$ there exists $\beta_\varepsilon < 1$ [independent of the sequence $\eta(\cdot)$] such that*

$$(46) \qquad \sup\{|\psi_k(u\eta(k))| : u \geq \varepsilon, k \geq 1\} < \beta.$$

Our next lemma gives the behavior of $\int_{0+}^{\pi} e^{-iu} f_r(\psi_s(u))\, du$ when $r, s \to \infty$ when $0 < \alpha < 1$.

LEMMA 3. *Assume $0 < \alpha < 1$ and $E(Z_1 \log Z_1) < \infty$. Let $\{\eta(r,s): r \geq 1, s \geq 1\}$ be a sequence of positive numbers such that $\inf_{r,s \geq 1} \eta(r,s) > 0$ and $\lim_{r,s \to \infty} \eta(r,s) = 1$. Then*

$$(47) \qquad \lim_{r,s \to \infty} p_1^{-r} \int_{0+}^{\pi} e^{-iu} f_r(\psi_s(u\eta(r,s)))\, du = \sum_{l \geq 0} (p_1 m)^{-l} I_l,$$

*where*

$$(48) \qquad I_l \equiv \int_{\pi m^{-1}}^{\pi} e^{-ium^{-l}} Q(\psi(u))\, du.$$

*Furthermore, there exists a universal constant $C$ such that*

$$(49) \qquad \sup_{r,s \geq 1} \left| p_1^{-r} \int_{0+}^{\pi} e^{-iu} f_r(\psi_s(u\eta(r,s)))\, du \right| \leq C.$$

PROOF. The finiteness of $\alpha$ renders $p_1 \neq 0$, and hence using the decomposition from Lemma 1 we have that

$$(50) \qquad \begin{aligned} p_1^{-r} &\int_{0+}^{\pi} e^{-iu} f_r(\psi_s(u\eta(r,s)))\, du \\ &= (p_1 m)^{-r} I^{(1)}(r,s) + \sum_{l=0}^{r} (p_1 m)^{-l} p_1^{-(r-l)} I_l^{(2)}(r,s), \end{aligned}$$

where

$$I^{(1)}(r,s) \equiv \int_0^{\pi m^{-1}} e^{-ium^{-r}} \psi_{s+r}(u\eta(r,s))\,du \tag{51}$$

and

$$I^{(2)}_l(r,s) \equiv \int_{\pi m^{-1}}^{\pi} e^{-ium^{-l}} f_{r-l}(\psi_{s+l}(u\eta(r,s)))\,du. \tag{52}$$

Note by uniform convergence on compacts of $\psi_r(u\eta(r,s))$, to $\psi(u)$ it follows that

$$\lim_{r,s\to\infty} I^{(1)}(r,s) = \int_0^{\pi m^{-1}} \psi(u)\,du. \tag{53}$$

Thus the first term on the right-hand side of (50) converges to 0 since $\alpha < 1$ is equivalent to $p_1 m > 1$. As for the second term, note first that

$$p_1^{-(r-l)} I^{(2)}_l(r,s) = \int_{\pi m^{-1}}^{\pi} e^{-ium^{-l}} Q_{r-l}(\psi_{s+l}(u\eta(r,s)))\,du. \tag{54}$$

Now

$$\left|\int_{\pi m^{-1}}^{\pi} e^{-ium^{-l}} Q_{r-l}(\psi_{s+l}(u\eta(r,s)))\,du\right|$$

$$\leq \int_{\pi m^{-1}}^{\pi} Q_{r-l}(|\psi_{s+l}(u\eta(r,s))|)\,du$$

$$\leq (\pi - \pi m^{-1}) Q_{r-l}(\beta),$$

where $\beta = \sup_{\pi m^{-1} \leq u \leq \pi, r, s \geq 1} |\psi_{s+l}(u\eta(r,s))| < 1$, by Lemma 2. Thus, by the monotonicity of $Q_r(\beta)$ in $r$,

$$\sup_{r,s} \sup_{l\leq r} \left|\int_{\pi m^{-1}}^{\pi} e^{-ium^{-l}} Q_{r-l}(\psi_{s+l}(u\eta(r,s)))\,du\right| \leq CQ(\beta) < \infty. \tag{55}$$

Thus, by the dominated convergence theorem, we have that

$$\lim_{r,s\to\infty} \sum_{l=0}^r (p_1 m)^{-l} p_1^{-(r-l)} I^{(2)}_l(r,s) \tag{56}$$
$$= \sum_{l\geq 0} (p_1 m)^{-l} \int_{\pi m^{-1}}^{\pi} e^{-ium^{-l}} Q(\psi(u))\,du,$$

thus concluding the proof of (47). Finally, note that by (50), (55) and the trivial estimate $|\psi(u)| \leq 1$,

$$\left|p_1^{-r} \int_{0+}^{\pi} e^{-iu} f_r(\psi_s(u\eta(r,s)))\,du\right| \leq (1 - p_1 m)^{-1} Q(\beta) + 1; \tag{57}$$

since the right-hand side of the above equation is independent of $r$ and $s$, (49) follows. □

Lemma 4 considers the case $\alpha = 1$. The proof of the lemma depends on the rate of convergence of $p_1^{-(r-l)} I_l(r,s)$ to $I_l$.

LEMMA 4. *Assume $\alpha = 1$ and $E(Z_1 \log Z_1) < \infty$. Let $\{\eta(r,s) : r \geq 1, s \geq 1\}$ be a sequence of positive numbers such that $\inf_{r,s \geq 1} \eta(r,s) > 0$ and $\lim_{r,s \to \infty} \eta(r,s) = 1$. Then*

$$(58) \quad \lim_{r,s \to \infty} r^{-1} p_1^{-r} \int_{0+}^{\pi} e^{-iu} f_r(\psi_s(u\eta(r,s))) \, du = \int_{\pi m^{-1}}^{\pi} Q(\psi(u)) \, du.$$

*Furthermore, there exists a universal constant $C$ such that*

$$(59) \quad \sup_{r,s \geq 1} \left| r^{-1} p_1^{-r} \int_{0+}^{\pi} e^{-iu} f_r(\psi_s(u\eta(r,s))) \, du \right| \leq C.$$

PROOF. From (50) and using $p_1 m = 1$, we have that

$$(60) \quad \begin{aligned} & r^{-1} p_1^{-r} \int_{0+}^{\pi} e^{-iu} f_r(\psi_s(u\eta(r,s))) \, du \\ & = r^{-1} I^{(1)}(r,s) + r^{-1} \sum_{l=0}^{r} p_1^{-(r-l)} I_l^{(2)}(r,s), \end{aligned}$$

where $p_1^{-(r-l)} I_l^{(2)}(r,s)$ is as in (54). From (53) it follows that $\lim_{r,s \to \infty} r^{-1} I^{(1)}(r,s) = 0$. Now for each fixed $l$, using (55) and the dominated convergence theorem,

$$(61) \quad \lim_{r,s \to \infty} \int_{\pi m^{-1}}^{\pi} e^{-ium^{-l}} Q_{r-l}(\psi_{s+l}(u\eta(r,s))) \, du = I_l,$$

where $I_l$ is as in (48). Furthermore, first by Lemma 2 and then the bounded convergence theorem,

$$(62) \quad \lim_{l \to \infty} I_l = \int_{\pi m^{-1}}^{\pi} Q(\psi(u)) \, du.$$

Thus,

$$(63) \quad r^{-1} \sum_{l=0}^{r} p_1^{-(r-l)} I_l^{(2)}(r,s) = r^{-1} \sum_{l=0}^{r} (p_1^{-(r-l)} I_l^{(2)}(r,s) - I_l) + r^{-1} \sum_{l=0}^{r} I_l.$$

By (62),

$$(64) \quad \lim_{r \to \infty} r^{-1} \sum_{l=0}^{r} I_l = \int_{\pi m^{-1}}^{\pi} Q(\psi(u)) \, du.$$

We will now show that $\lim_{r,s \to \infty} r^{-1} \sum_{l=0}^{r} (p_1^{-(r-l)} I_l^{(1)}(r,s) - I_l) = 0$. To this end,

$$(65) \quad r^{-1} \sum_{l=0}^{r} (p_1^{-(r-l)} I_l^{(2)}(r,s) - I_l) = T_1(r,s) + T_2(r,s),$$

where

$$T_1(r,s) \equiv r^{-1} \sum_{l=0}^{r} I_{l,1}^{(2)}(r,s),$$

$$T_2(r,s) = r^{-1} \sum_{l=0}^{r} I_{l,2}^{(2)}(r,s),$$

(66)    $$I_{l,1}^{(2)}(r,s) = \int_{\pi m^{-1}}^{\pi} e^{-ium^{-l}} (Q_{r-l}(\psi_{s+l}(u\eta(r,s))) - Q_{r-l}(\psi(u)))\, du$$

and

(67)    $$I_{l,2}^{(2)}(r,s) \equiv \int_{\pi m^{-1}}^{\pi} e^{-ium^{-l}} (Q_{r-l}(\psi(u)) - Q(\psi(u)))\, du.$$

We first show $\lim_{r,s\to\infty} T_2(r,s) = 0$. By using the change of variable in the index of summation, one has

(68)    $$|T_2(r,s)| \le (\pi - \pi m^{-1}) r^{-1} \sum_{l=0}^{r} \sup_{\pi m^{-1} \le u \le \pi} |Q_{r-l}(\psi(u)) - Q(\psi(u))|;$$

now since $\sup\{|\psi(u)| : \pi m^{-1} \le u \le \pi\} < 1$, letting $r,s \to \infty$ and using the uniform convergence of $Q_r(\cdot)$ to $Q(\cdot)$ in the interior of the unit disk one has that $\lim_{r,s\to\infty} T_2(r,s) = 0$. In fact, it also follows that there exists a constant $C$ such that

(69)    $$\sup_{r,s \ge 1} |T_2(r,s)| \le C.$$

We now move on to establish that $\lim_{r,s\to\infty} T_1(r,s) = 0$. Note that, by monotonicity of $Q'_r(s)$ in $r$,

$$|T_1(r,s)| \le r^{-1} \sum_{l=0}^{r} \int_{\pi m^{-1}}^{\pi} |Q_{r-l}(\psi_{s+l}(u\eta(r,s))) - Q_{r-l}(\psi(u))|\, du$$

$$\le r^{-1} Q'(\beta) \sum_{l=0}^{r} \int_{\pi m^{-1}}^{\pi} |\psi_{s+l}(u) - \psi(u)|\, du,$$

where $\beta \equiv \sup\{|\delta\psi_s(u\eta(r,s)) + (1-\delta)\psi(u)| : 0 \le \delta \le 1, \pi m^{-1} \le u \le \pi, r, s \ge 1\} < 1$ by Lemma 2. Thus, $T_1(r,s)$ converges to 0 as $r,s \to \infty$ due to the uniform convergence, in the interval $[\pi m^{-1}, \pi]$, of $\psi_r(u)$ to $\psi(u)$ as $r \to \infty$. The above calculation also yields

(70)    $$\sup_{r,s \ge 1} |T_1(r,s)| \le C$$

for some universal constant $C$. Thus (58) follows by letting $r, s \to \infty$ in (60) and using (64). Further (59) also follows since by (69) and (70),

$$\sup_{r,s \geq 1} \left| r^{-1} \sum_{l=0}^{r} p_1^{-(r-l)} I_l^{(2)}(r,s) \right| \leq C \tag{71}$$

and

$$\sup_{r,s \geq 1} |r^{-1} I^{(1)}(r,s)| \leq 2. \tag{72}$$

$\square$

Our next lemma gives the behavior of $\int_{0+}^{\pi} e^{-iu} f_r(\psi_s(u)) \, du$ when $r, s \to \infty$ when $\alpha > 1$.

LEMMA 5. *Assume $\alpha > 1$ and $E(Z_1 \log Z_1) < \infty$. Let $\{\eta(r,s): r \geq 1, s \geq 1\}$ be a sequence of positive numbers such that $\inf_{r,s \geq 1} \eta(r,s) > 0$ and $\lim_{r,s \to \infty} \eta(r,s) = 1$. Then*

$$\lim_{r,s \to \infty} m^r \int_{0+}^{\pi} e^{-iu} f_r(\psi_s(u\eta(r,s))) \, du = K, \tag{73}$$

*where*

$$K = \begin{cases} \sum_{l \geq 0} (p_1 m)^l \int_{\pi m^{-1}}^{\pi} Q_l(\psi(u)) \, du + \int_{0}^{\pi m^{-1}} \psi(u) \, du, & \text{if } \alpha < \infty, \\ \int_{0}^{\pi m^{-1}} \psi(u) \, du, & \text{if } \alpha = \infty. \end{cases}$$

*Furthermore, there exists a universal constant $C$ such that*

$$\sup_{r,s \geq 1} \left| m^r \int_{0+}^{\pi} e^{-iu} f_r(\psi_s(u\eta(r,s))) \, du \right| \leq C. \tag{74}$$

PROOF. When $\alpha < \infty$, multiplying (50) by $(p_1 m)^r$ one gets

$$\left| m^r \int_{0+}^{\pi} e^{-iu} f_r(\psi_s(u\eta(r,s))) \, du \right| \\ = I^{(1)}(r,s) + \sum_{l=0}^{r} (p_1 m)^{r-l} p_1^{-(r-l)} I_l^{(2)}(r,s), \tag{75}$$

where $p_1^{-(r-l)} I_l^{(2)}(r,s)$ and $I^{(1)}(r,s)$ are as defined in (52) and (51), respectively. Now,

$$T_1(r,s) = \sum_{l=0}^{r} (p_1 m)^{r-l} p_1^{-(r-l)} I_l^{(1)}(r,s)$$

$$= \sum_{l \geq 0} (p_1 m)^l \left( \int_{\pi m^{-1}}^{\pi} e^{-ium^{r-l}} Q_l(\psi_{s+r-l}(u\eta(r,s))) \, du \right) I_{l \leq r}.$$

Thus, by the dominated convergence theorem [using (55) and $p_1 m < 1$],

$$
\text{(76)} \quad \lim_{r,s\to\infty} T_1(r,s) = \sum_{l\geq 0}(p_1 m)^l \int_{\pi m^{-1}}^{\pi} Q_l(\psi(u))\,du.
$$

Finally (73) follows by letting $r, s \to \infty$ in (75) and using (53). We now consider the case $\alpha = \infty$. Using (37) and multiplying by $m^r$ one gets

$$
\text{(77)} \quad m^r \int_0^{\pi} e^{-iu} f_r(\psi_s(u\eta(r,s)))\,du = I^{(1)}(r,s) + \sum_{l=0}^{r} m^{r-l} I_l^{(2)}(r,s),
$$

where

$$
\text{(78)} \quad I^{(1)}(r,s) \equiv \int_0^{\pi m^{-1}} e^{-ium^{-r}} \psi_{s+r}(u\eta(r,s))\,du
$$

and

$$
\text{(79)} \quad I_l^{(2)}(r,s) \equiv \int_{\pi m^{-1}}^{\pi} e^{-ium^{-l}} f_{r-l}(\psi_{s+l}(u\eta(r,s)))\,du.
$$

Using (53), it follows that $\lim_{r,s\to\infty} I^{(1)}(r,s) = \int_0^{\pi m^{-1}} \psi(u)\,du$. We now have to deal with $I_l^{(2)}(r,s)$; since $\alpha = \infty$, let $j_0 = \inf\{j \geq 2 : p_j > 0\}$. Then for any $\beta < 1$, we have

$$
\text{(80)} \quad f_n(\beta) \leq C\beta^{(j_0)^n};
$$

the above estimate can easily be proved from Proposition 3 of [1] or from Theorem 24.1 of [14]. Thus, it follows that there exists a universal constant $C$ such that

$$
\text{(81)} \quad |I_l^{(2)}(r,s)| \leq C\beta^{(j_0)^{r-l}},
$$

where $\beta = \sup\{|\psi_k(u)| : \eta_0 \pi m^{-1} \leq u \leq \pi, k \geq 1\} < 1$ and $\eta_0 = \inf\{\eta(r,s) : r, s \geq 1\}$ (which is positive by assumption). Now, using the estimate (81) and the dominated convergence theorem, it follows that

$$
\text{(82)} \quad \lim_{r,s\to\infty} \sum_{l=0}^{r} m^{r-l} I_l^{(2)}(r,s) = 0. \qquad \square
$$

The next lemma provides a uniform estimate on $\int_\pi^{\pi m} e^{-iut} Q_r(\psi_s(u))\,du$.

LEMMA 6. *Let $\{\eta(r,s) : r \geq 1,\ s \geq 1\}$ be a sequence of positive numbers such that $\inf_{r,s\geq 1} \eta(r,s) > 0$ and $\lim_{r,s\to\infty} \eta(r,s) = 1$. There exists a universal constant $C$ such that*

$$
\text{(83)} \quad \sup_{r,s\geq 1}\left|\int_\pi^{\pi m} e^{-iut} Q_r(\psi_s(u\eta(r,s)))\,du\right| \leq \frac{C}{|t|}.
$$

PROOF. The proof uses the following estimate whose proof can be constructed using methods similar to the one in [31], page 45, (4.2). Dubuc ([13], page 481) and Dubuc and Senata ([15], page 494) have used this estimate extensively in their work. For any $h \in L_\infty[a,b]$,

$$\text{(84)} \qquad \left| \int_a^b e^{iut} h(u)\, du \right| \leq \frac{\|h\|_\infty \pi}{|t|} + \frac{(b-a)}{2} w_h\left(\frac{\pi}{|t|}\right),$$

where $w_h(\delta)$ is the modulus of continuity of $h$ obtained using intervals of length $\delta$. Letting $h = \tilde{Q}_{r,s} \equiv Q_r(\psi_s(u\eta(r,s))$, we have

$$\left| \int_\pi^{\pi m} e^{-iut} Q_r(\psi_s(u\eta(r,s)))\, du \right|$$

$$\leq \frac{\|\tilde{Q}_{r,s}\|_\infty \pi}{|t|} + \frac{\pi(m-1)}{2} w_{\tilde{Q}_{r,s}}\left(\frac{\pi}{|t|}\right)$$

$$\leq \frac{Q(\beta)\pi}{|t|} + \frac{\pi(m-1)}{2} w_{\tilde{Q}_{r,s}}\left(\frac{\pi}{|t|}\right),$$

where $\beta = \sup\{|\psi_k(u)| : \pi\eta_0 \leq u \leq m\pi, k \geq 1\} < 1$ (by yet another application of Lemma 2) and $\eta_0 = \inf_{r,s \geq 1} \eta(r,s)$ which is positive by assumption. Also by the mean value theorem and the estimate

$$\text{(85)} \qquad \sup\{|\psi_s(u_1) - \psi_s(u_2)| : \eta_0\pi \leq u_1 \leq u_2 \leq m\pi, \\ |u_1 - u_2| \leq \pi|t|^{-1}\} \leq C|t|^{-1},$$

it follows that there exists a universal constant $C$ such that

$$\text{(86)} \qquad w_{\tilde{Q}_{r,s}}\left(\frac{\pi}{|t|}\right) \leq \frac{C}{|t|}.$$

Finally the estimate (83) follows from the uniform continuity of $\psi_s(u)$ and the uniform estimate $\sup\{|\psi'_s(u)| : \eta_0\pi \leq u_1 \leq u_2 \leq m\pi, s \geq 1\} \leq C$. □

Our next lemma considers the the behavior of $\int_\pi^{\pi m^s} e^{-iu} f_r(\psi_s(u))\, du$.

LEMMA 7. *Assume $E(Z_1 \log Z_1) < \infty$. Let $\{\eta(r,s) : r \geq 1, s \geq 1\}$ be a sequence of positive numbers such that $\inf_{r,s \geq 1} \eta(r,s) > 0$ and $\lim_{r,s \to \infty} \eta(r,s) = 1$. Then for any $0 < \alpha < 1$,*

$$\text{(87)} \quad \lim_{r,s \to \infty} p_1^{-r} \int_\pi^{\pi m^s} e^{-iu} f_r(\psi_s(u))\, du = \sum_{l \geq 0} (p_1 m)^l \int_\pi^{\pi m} e^{-ium^l} Q(\psi(u))\, du.$$

*Furthermore, there exists a universal constant $C$ such that*

$$\text{(88)} \qquad \sup_{r,s \geq 1} p_1^{-r} \left| \int_\pi^{\pi m^s} e^{-iu} f_r(\psi_s(u))\, du \right| \leq C.$$

PROOF. We begin with the decomposition similar to the one in Lemma 1, that is,

$$\int_\pi^{\pi m^s} e^{-iu} f_r(\psi_s(u\eta(r,s)))\, du = \sum_{l=0}^{s-1} \int_{\pi m^l}^{\pi m^{l+1}} e^{-iu} f_r(\psi_s(u\eta(r,s)))\, du$$

$$= \sum_{l=0}^{s-1} m^{-l} \int_{\pi m^l}^{\pi m^{l-1}} e^{-ium^l} f_{r+l}(\psi_{s-l}(u\eta(r,s)))\, du,$$

where the last identity follows from the change of variable $u \longmapsto um^r$ and the identity $\psi_n(sm^l) = f_l(\psi_{n-l}(s))$. Thus,

$$(89) \quad p_1^{-r} \int_\pi^{\pi m^s} e^{-iu} f_r(\psi_s(u))\, du = \sum_{l=0}^{s-1} (p_1 m)^l \int_\pi^{\pi m} e^{-ium^l} Q_{r+l}(\psi_{s-l}(u))\, du.$$

Now (88) follows from Lemma 6. Finally, Lemma 6 and the dominated convergence theorem yield (87). □

LEMMA 8. *Assume* $E(Z_1 \log Z_1) < \infty$ *and* $\alpha = 1$. *Let* $\{\eta(r,s): r \geq 1, s \geq 1\}$ *be a sequence of positive numbers such that* $\inf_{r,s \geq 1} \eta(r,s) > 0$ *and* $\lim_{r,s \to \infty} \eta(r,s) = 1$. *Then,*

$$(90) \quad \lim_{r,s \to \infty} s^{-1} p_1^{-r} \int_\pi^{\pi m^s} e^{-iu} f_r(\psi_s(u\eta(r,s)))\, du = 0.$$

*Furthermore, there exists a universal constant* $C$ *such that*

$$(91) \quad \sup_{r,s \geq 1} s^{-1} p_1^{-r} \left| \int_\pi^{\pi m^s} e^{-iu} f_r(\psi_s(u\eta(r,s)))\, du \right| \leq C.$$

PROOF. Using (89) and (72),

$$(92) \quad \begin{aligned} & s^{-1} p_1^{-r} \int_\pi^{\pi m^s} e^{-iu} f_r(\psi_s(u\eta(r,s)))\, du \\ &= \frac{1}{s} \sum_{l=0}^{s-1} \int_\pi^{\pi m} e^{-ium^l} Q_{r+l}(\psi_{s-l}(u\eta(r,s)))\, du \end{aligned}$$

and the (90) follows by Lemma 6 and taking the limits as $s \to \infty$. Equation (91) follows from Lemma 6. □

Our next lemma considers the case $\alpha > 1$.

LEMMA 9. *Assume* $E(Z_1 \log Z_1) < \infty$. *Let* $\{\eta(r,s): r \geq 1, s \geq 1\}$ *be a sequence of positive numbers such that* $\inf_{r,s \geq 1} \eta(r,s) > 0$ *and* $\lim_{r,s \to \infty} \eta(r,s) = 1$. *Then for any* $\alpha > 1$,

$$(93) \quad \lim_{r,s \to \infty} m^r \int_\pi^{\pi m^s} e^{-iu} f_r(\psi_s(u))\, du = 0.$$

*Furthermore, there exists a universal constant $C$ such that*

$$\text{(94)} \qquad \sup_{r,s\geq 1}\left|m^r \int_\pi^{\pi m^s} e^{-iu} f_r(\psi_s(u))\, du\right| \leq C.$$

PROOF. When $p_1 > 0$ multiplying (89) by $(mp_1)^r$ and using Lemma 6, (93) and (94) follows. When $\alpha = \infty$, then using the decomposition

$$\text{(95)} \qquad \begin{aligned} & m^r \int_\pi^{\pi m^s} e^{-iu} f_r(\psi_s(u\eta(r,s)))\, du \\ &= m^r \sum_{l=0}^{s-1} m^l \int_\pi^{\pi m} e^{-ium^l} f_{r+l}(\psi_{s-l}(u\eta(r,s)))\, du \end{aligned}$$

and the estimate

$$\text{(96)} \qquad \left|\int_\pi^{\pi m} e^{-ium^l} f_{r+l}(\psi_{s-l}(u\eta(r,s)))\, du\right| \leq C\beta^{(j_0)^{r+l}}$$

[where $\beta = \sup\{|\psi_k(u)| : \pi\eta_0 \leq u \leq m\pi, k \geq 1\}$ ($< 1$ by Lemma 2), $j_0 = \inf\{j \geq 2 : p_j > 0\}$ and $C$ is a universal constant] one can complete the proof of the lemma. $\square$

We are now ready to give a proof of Theorem 1.

PROOF OF THEOREM 1. Using the inversion theorem and change of variables it follows that

$$\text{(97)} \qquad (2\pi)v_n P(Z_n = v_n) = J_n^1 + J_n^2,$$

where

$$\text{(98)} \qquad J_n^1 = \int_{-\pi v_n}^{0-} e^{-iu} f_{k_n}(\psi_{a_n}(uv_n^{-1} m^{a_n}))\, du$$

and

$$\text{(99)} \qquad J_n^2 = \int_{0+}^{\pi v_n} e^{-iu} f_{k_n}\left(\psi_{a_n}\left(u\frac{m^{a_n}}{v_n}\right)\right) du,$$

where $k_n = [\tilde{k}_n]$, $\tilde{k}_n$ is such that $m^{n-\tilde{k}_n} = v_n$, and $a_n = n - k_n$. We will only deal with $J_n^2$ since $J_n^1$ can be handled similarly. First, decompose

$$\text{(100)} \qquad J_n^2 = \left(\int_{0+}^\pi + \int_\pi^{\pi v_n}\right) e^{-iu} f_{k_n}\left(\psi_{a_n}\left(u\frac{m^{a_n}}{v_n}\right)\right) du.$$

*Case $\alpha < 1$.* When $\alpha < 1$, we have by Lemma 3.3, when $\frac{m^{a_n}}{v_n}$ converges to 1,

$$\lim_{n\to\infty} p_1^{k_n} \int_{0+}^\pi e^{-iu} f_{k_n}\left(\psi_{a_n}\left(u\frac{m^{a_n}}{v_n}\right)\right) du$$

$$(101) \qquad = \sum_{l \geq 0} (p_1 m)^{-l} I_l < \infty$$

$$(102) \qquad = \int_0^\pi e^{-iu} Q(\psi(u)) \, du,$$

where $I_l$ is defined as in (48); the last equality follows from the change of variable $um^{-l} \longmapsto u$, and using $\psi(um^l) = f_l(\psi(u))$, and the functional equation $Q(f_l(\psi(u))) = p_1^l Q(\psi(u))$. Now using Lemma 6 and similar arguments used in establishing (102), it follows that

$$(103) \qquad \lim_{n \to \infty} p_1^{k_n} \int_\pi^{\pi v_n} e^{-iu} f_{k_n}\left(\psi_{a_n}\left(u \frac{m^{a_n}}{v_n}\right)\right) du$$

$$= \sum_{l \geq 0} (p_1 m)^l \int_\pi^{m\pi} e^{-ium^l} Q(\psi(u)) \, du$$

$$(104) \qquad = \int_\pi^\infty e^{-iu} Q(\psi(u)) \, du < \infty.$$

Thus, using $p_1^{k_n} = p_1^n v_n^\alpha$, the definition of $Q(\cdot)$, and the inversion theorem, we have that

$$(105) \qquad \lim_{n \to \infty} \frac{P(Z_n = v_n)}{p_1^n v_n^{\alpha - 1}} = \frac{1}{2\pi} \int_{-\infty}^\infty e^{-iu} Q(\psi(u)) \, du$$

$$(106) \qquad = \sum_{j \geq 1} q_j w^{\star j}(1) < \infty,$$

thus proving the theorem when $\frac{m^{a_n}}{v_n}$ converges to 1. Now when $\frac{m^{a_n}}{v_n}$ does not converge, it oscillates between 1 and $m$. Thus we have for $1 < c < m$

$$\inf_{1 \leq c \leq m} p_1^{-k_n}(J_n^1(c) + J_n^2(c)) \leq \frac{P(Z_n = v_n)}{p_1^n v_n^{\alpha-1}} \leq \sup_{1 \leq c \leq m} p_1^{-k_n}(J_n^1(c) + J_n^2(c)),$$

where

$$(107) \qquad J_n^1(c) = \int_{-\pi v_n}^{0-} e^{-iu} f_{k_n}(\psi_{a_n}(uc)) \, du$$

and

$$(108) \qquad J_n^2(c) = \int_{0+}^{\pi v_n} e^{-iu} f_{k_n}(\psi_{a_n}(uc)) \, du.$$

Using (101)–(106) it follows that for every $1 \leq c \leq m$,

$$(109) \qquad \lim_{n \to \infty} p_1^{-k_n}(J_n^1(c) + J_n^2(c)) = \frac{1}{c} \sum_{j \geq 1} q_j w^{\star j}(c) < \infty.$$

Thus, with $C_1 = \inf_{1 \leq c \leq m} \frac{1}{c} \sum_{j \geq 1} q_j w^{\star j}(c)$ and $C_2 = \sup_{1 \leq c \leq m} \frac{1}{c} \sum_{j \geq 1} q_j w^{\star j}(c)$, we have

$$(110) \quad C_1 \leq \liminf_{n \to \infty} \frac{P(Z_n = v_n)}{p_1^n v_n^{\alpha-1}} \leq \limsup_{n \to \infty} \frac{P(Z_n = v_n)}{p_1^n v_n^{\alpha-1}} \leq C_2,$$

thus concluding the proof for Case $\alpha < 1$.

*Case $\alpha = 1$.* When $\alpha = 1$, we have by Lemma 4, if $\frac{m^{a_n}}{v_n} \to 1$, that

$$\lim_{n \to \infty} p_1^{k_n} k_n^{-1} \int_{0+}^{\pi} e^{-iu} f_{k_n}\left(\psi_{a_n}\left(u \frac{m^{a_n}}{v_n}\right)\right) du = \int_{\pi m^{-1}}^{\pi} Q(\psi(u)) \, du,$$

where $|\int_{\pi m^{-1}}^{\pi} Q(\psi(u)) \, du| < \infty$. Furthermore, by Lemma 8,

$$\lim_{n \to \infty} k_n^{-1} p_1^{k_n} \int_{\pi}^{\pi v_n} e^{-iu} f_{k_n}\left(\psi_{a_n}\left(u \frac{m^{a_n}}{v_n}\right)\right) du = 0,$$

since $\frac{a_n}{k_n} \leq 1$. Thus, using $p_1^{k_n} = p_1^n v_n$ and that $\frac{m^{a_n}}{v_n} \to 1$, we have

$$\lim_{n \to \infty} \frac{P(Z_n = v_n)}{p_1^n k_n} = \int_{-\pi}^{-\pi m^{-1}} Q(\psi(u)) \, du + \int_{\pi m^{-1}}^{\pi} Q(\psi(u)) \, du$$

$$= \int_{\pi m^{-1}}^{\pi} (Q(\psi(u)) - Q(\psi(-u))) \, du.$$

The finiteness of the limit follows from Lemmas 4 and 7. To complete the proof we need to establish the positivity of the limit. To this end, first note that there exists an interval $(a,b) \in (\pi m^{-1}, \pi)$ such that $\int_a^b (Q(\psi(u)) - Q(\psi(-u))) \, du \neq 0$; for, if not, that would imply $\psi^j(t) = \psi^j(-t)$ for all $t \in (\pi m^{-1}, \pi)$ and all $j \geq 1$, which is impossible. Let $(a,b)$ be the largest such interval. Then $\int_{\pi m^{-1}}^{\pi} (Q(\psi(u)) - Q(\psi(-u))) \, du = \int_a^b (Q(\psi(u)) - Q(\psi(-u))) \, du \neq 0$. Now when, $\frac{m^{a_n}}{v_n}$ does not have a limit as $n \to \infty$, we proceed exactly as in Case $\alpha < 1$.

*Case $\alpha > 1$.* When $\alpha > 1$, we have by Lemma 5, $\frac{m^{a_n}}{v_n} \to 1$, that

$$\lim_{n \to \infty} m^{k_n} \int_{0+}^{\pi} e^{-iu} f_{k_n}\left(\psi_{a_n}\left(u \frac{m^{a_n}}{v_n}\right)\right) du = K,$$

where $|K| < \infty$ and is given in Lemma 5. Furthermore, by Lemma 7,

$$\lim_{n \to \infty} m^{k_n} \int_{\pi}^{\pi v_n} e^{-iu} f_{k_n}\left(\psi_{a_n}\left(u \frac{m^{a_n}}{v_n}\right)\right) du = 0.$$

When $v_n^{-1} m^{a_n}$ does not converge to 1 the result is established arguing as in Cases $\alpha < 1$ and $\alpha = 1$. □

PROOF OF THEOREM 2. The lower bound follows from Theorem 1. Thus it is enough to establish the upper bound. To this end, we first note that

$$(111) \quad \sup_{j \geq v_n} P(Z_n = j) = \sup_{r \geq 1} \sup_{j \leq v_n} P(Z_n = rv_n + j).$$

We will obtain estimates of $P(Z_n = rv_n + j)$. Using the inversion formula and change of variables we have

$$(2\pi)(rv_n + j)P(Z_n = rv_n + j) = J_n^1 + J_n^2, \tag{112}$$

where

$$J_n^1 = \int_{-\pi(rv_n+j)}^{0-} e^{-iu} f_{k_n}(\psi_{a_n}(u\phi(n,j)))\, du, \tag{113}$$

$$J_n^2 = \int_{0+}^{\pi(rv_n+j)} e^{-iu} f_{k_n}(\psi_{a_n}(u\phi(n,j)))\, du, \tag{114}$$

$$\phi(n,j) = \left(\frac{m^{a_n}}{rv_n}\right)\left(\frac{rv_n}{rv_n + j}\right), \tag{115}$$

where $m^{\tilde{k}_n} = rv_n$ and $k_n = [\tilde{k}_n]$ and $a_n = n - k_n$. We first decompose

$$J_n^2 = J_{n,A}^2 + J_{n,B}^2 + J_{n,C}^2, \tag{116}$$

where

$$J_{n,A}^2 = \int_{0+}^{\pi} e^{-iu} f_{k_n}(\psi_{a_n}(u\phi(n,j)))\, du, \tag{117}$$

$$J_{n,B}^2 = \int_{\pi}^{\pi(rv_n)} e^{-iu} f_{k_n}(\psi_{a_n}(u\phi(n,j)))\, du \tag{118}$$

and

$$J_{n,C}^2 = \int_{\pi rv_n}^{\pi(rv_n+j)} e^{-iu} f_{k_n}(\psi_{a_n}(u\phi(n,j)))\, du. \tag{119}$$

A similar decomposition also holds for $J_n^1$.

*Case $\alpha < 1$.* Using the uniform estimate from Lemmas 4 and 7 it follows that there exists a universal constant $C_1$ (independent of $r, n, j$) such that

$$(J_{n,A}^2 + J_{n,B}^2) < C_1 p_1^{k_n}. \tag{120}$$

Now, it can be shown that

$$J_{n,C}^2 = p_1^{n-1} m^{a_n - 1} \int_{m\pi}^{m\pi + jm^{-(a_n+1)}} e^{-ium^{a_n+1}} Q_{n-1}(\psi_1(u\phi(n,j)))\, du. \tag{121}$$

Thus, using the by now standard arguments, it follows that there exists a universal constant $C$ such that

$$|J_{n,C}^2| \leq C p_1^{n-1} m^{a_n - 1} \leq p_1^{k_n}. \tag{122}$$

Thus

$$|J_n^2| \leq p_1^{k_n}, \tag{123}$$

and a similar estimate holds for $J_n^1$. Thus, it follows that (using the definition of $k_n$)

$$(124) \qquad P(Z_n = rv_n + j) \leq C_2 p_1^n (rv_n + j)^{\alpha-1},$$

where $C_2$ is a universal constant.

*Case $\alpha = 1$.* Using the uniform estimate from Lemmas 5, 8 and an argument similar to Case $\alpha < 1$, it follows that there exist a universal constant $C_3$ (independent of $r, n, j$) such that

$$(125) \qquad (2\pi)(rv_n + j)P(Z_n = rv_n + j) < C_3 p_1^{k_n} k_n.$$

Now using the definition of $k_n$ it follows that

$$(126) \qquad P(Z_n = rv_n + j) \leq C_4 p_1^n k_n,$$

where $C_4$ is a universal constant.

*Case $\alpha > 1$.* Using the uniform estimate from Lemmas 6 and 9, and an argument similar to Case $\alpha < 1$, it follows that there exists a universal constant $C_5$ (independent of $r, n, j$) such that

$$(127) \qquad (2\pi)(rv_n + j)P(Z_n = rv_n + j) < C_5 m^{k_n}.$$

Now using the definition of $k_n$ it follows that

$$(128) \qquad P(Z_n = rv_n + j) \leq C_6 m^{-n},$$

where $C_6$ is a universal constant. □

PROOF OF PROPOSITION 1. Follows from the proof below with $v_{n-k} = 1$. □

PROOF OF PROPOSITION 2. Using independence of $Y_n$ and $N_n$, the definition of conditional probability and the definition of the random variable $N_n^1(k)$, it follows that

$$(129) \qquad P(Y_{N_n} \in A | N_{n-k} \geq v_{n-k}) = \sum_{l \geq 0} P(Y_l \in A) P(N_n^1(k) = l).$$

Now, since $\{Y_n : n \geq 1\}$ satisfies an LDP with speed $r(n)$, there exist constants $C$ and $0 < \delta < 1$ such that

$$P(Y_{N_n} \in A | N_{n-k} \geq v_{n-k}) \leq C \sum_{l \geq 0} e^{r(l)(I_Y(\bar{A}) - \delta)} P(N_n^1(k) = l)$$

$$= CE(\exp(-\phi r(N_n^1(k)))),$$

where $\phi = I_Y(\bar{A}) - \delta$ and $\bar{A}$ is the closure of the set $A$. The result follows by taking the logarithm, using (17) and letting $\delta \searrow 0$. The lower bound is proved similarly. The fact that $K(I_Y(\cdot))$ is lower semicontinuous follows from the assumed continuity of $K(\cdot)$ and the lower semicontinuity of $I_Y(\cdot)$.
□

We now move to prove the second main result of this paper, namely Theorem 3.

PROOF OF THEOREM 3. Using Proposition 2, it is enough to evaluate for $\theta > 0$,

$$\lim_{n \to \infty} \frac{1}{v_{n-k}} \log E(\exp(-\theta N_n^1(k))) \equiv K(\theta) \tag{130}$$

and check that it is continuous, where $\{N_n^1(k) : n \geq k\}$ is a sequence of random variables with distribution

$$P(N_n^1(k) = l) = P(Z_n = l | Z_{n-k} \geq v_{n-k}). \tag{131}$$

Using the definition of conditional expectations one can show that

$$\begin{aligned}
E(\exp(-\theta N_n^1(k))) &\\
&= E(\exp(-\theta Z_n) | Z_{n-k} \geq v_{n-k}) \\
&= \sum_{j \geq v_{n-k}} (f_k(e^{-\theta}))^j \frac{P(Z_{n-k} = j)}{P(Z_{n-k} \geq v_{n-k})} \\
&= \frac{(f_k(e^{-\theta}))^{v_{n-k}}}{P(Z_{n-k} \geq v_{n-k})} \sum_{j \geq 0} (f_k(e^{-\theta}))^j P(Z_{n-k} = j + v_{n-k}) \\
&= \frac{(f_k(e^{-\theta}))^{v_{n-k}} A_{n-k}}{P(Z_{n-k} \geq v_{n-k})} \sum_{j \geq 0} (f_k(e^{-\theta}))^j \frac{P(Z_{n-k} = j + v_{n-k})}{A_{n-k}},
\end{aligned}$$

where $A_n$ is as given in Theorem 1. We now use Theorem 2 to establish the following.

CLAIM 1. *For every $0 < \alpha < \infty$,*

$$\lim_{n \to \infty} \frac{1}{v_{n-k}} \log \sum_{j \geq 0} (f_k(e^{-\theta}))^j \frac{P(Z_{n-k} = j + v_{n-k})}{A_{n-k}} = 0. \tag{132}$$

PROOF. Case $0 < \alpha < 1$.

$$\begin{aligned}
\sum_{j \geq 0} (f_k(e^{-\theta}))^j &\frac{P(Z_{n-k} = j + v_{n-k})}{A_{n-k}} \\
&= \sum_{j \geq 0} (f_k(e^{-\theta}))^j (j + v_{n-k})^{\alpha - 1} \frac{P(Z_{n-k} = j + v_{n-k})}{A_{n-k}(j + v_{n-k})^{\alpha - 1}}
\end{aligned}$$

$$\leq (v_{n-k})^{1-\alpha} \sum_{j\geq 0}(f_k(e^{-\theta}))^j \frac{\sup_{j\geq v_{n-k}} j^{\alpha-1}P(Z_{n-k}=j)}{A_{n-k}}$$

$$\leq C(v_{n-k})^{1-\alpha} \sum_{j\geq 0}(f_k(e^{-\theta}))^j \qquad \text{by Theorem 2.}$$

Claim 1 follows by the finiteness of $\sum_{j\geq 0}(f_k(e^{-\theta}))^j$ for $\theta > 0$. The other two cases follow a similar pattern of proof. □

Thus,

$$(133) \qquad \lim_{n\to\infty} \frac{1}{v_{n-k}} \log E(\exp(-\theta N_n^1(k))) = \log f_k(e^{-\theta}) - bB \equiv K(\theta),$$

where we have used that $\lim_{n\to\infty} \frac{1}{v_n} \log A_n = bB$ and $\lim_{n\to\infty} \frac{1}{v_n} \log P(Z_n \geq v_n) = 0$ since $v_n = \mathbf{O}(m^n)$ as $n \to \infty$. By continuity of $f_k(\cdot)$ it follows that $K(\theta)$ is continuous and hence the rate function is, by Proposition 2, $-K(I(x)) = -\log f_k(e^{-I(x)}) + bB \equiv \tilde{I}(x)$.

Now for any $L < \infty$,

$$\begin{aligned} V_L &\equiv \{x|\tilde{I}_k(x) \leq L\} \\ &= \{x|-\log f_k(e^{-I(x)}) + bB \leq L\} \\ &= \{x|I(x) \leq -\log g_k(e^{-(L-bB)})\}, \end{aligned}$$

where $g_k(\cdot)$ is the $k$th iterate of $g(\cdot)$, the functional inverse of $f(\cdot)$. Thus if $I(\cdot)$ is a good rate function, $\tilde{I}(\cdot)$ is also a good rate function. □

PROOF OF COROLLARIES 1 AND 2. Since $R_n = \frac{1}{Z_n}\sum_{i=1}^{Z_n} \xi_{n,i}$, the LDP for $\mu_{n,k}^{(1)}$ follows by taking $Y_n = \frac{1}{n}\sum_{i=1}^n \xi_{n,i}$. The LDP for $\mu_{n,k}^{(2)}$ follows by noting that

$$(134) \qquad \frac{W}{W_n} = \frac{1}{Z_n} \sum_{i=1}^{Z_n} W^{(j)},$$

where $W^{(j)}$'s are i.i.d. random variables distributed as $W$ and taking $Y_n = \frac{1}{n}\sum_{i=1}^n W^{(j)}$. □

PROOF OF COROLLARY 3. The proof follows by conditioning and using the estimates as in [11], page 176. □

PROOF OF THEOREM 4. By Proposition 2 it is enough to show that for $\theta > 0$ and $b > 0$,

$$\lim_{n\to\infty} \frac{1}{v_{n-k}} \log E(e^{-\theta N_n^1(k)}) \equiv bG(f_k(\theta))$$

and that $G(f_k(\theta))$ is continuous in $\theta$, where $N_n^1(k)$ is a random variable with distribution

$$P(N_n^1(k) = j) = P(Z_n = j | Z_{n-k} \geq v_{n-k}).$$

Now for $\theta > 0$, setting $\delta_n = P(Z_n \geq v_n)$, we have

$$E(e^{\theta N_n^1(k)}) = \sum_{j \geq 0} e^{-\theta j} P(Z_n = j | Z_{n-k} \geq v_{n-k})$$

$$= (\delta_{n-k})^{-1} \sum_{j \geq v_{n-k}} (f_k(e^{-\theta}))^r P(Z_{n-k} = r);$$

(135)
$$= (\delta_{n-k})^{-1} f_{n-k}(f_k(e^{-\theta})) - \sum_{j=j_0^{n-k}}^{v_{n-k}-1} (f_k(e^{-\theta}))^r P(Z_{n-k} = r)$$

$$= (f_k(e^{-\theta}))^{j_0^{n-k}} (\delta_{n-k})^{-1}$$
$$\times (f_{n-k}(f_k(e^{-\theta}))(f_k(e^{-\theta}))^{-j_0^{n-k}} - A(n,k)),$$

where $A(n,k) = \sum_{j=0}^{v_{n-k}-j_0^{n-k}-1} (f_k(e^{-\theta}))^r P(Z_{n-k} = r + j_0^{n-k})$.

CLAIM 2.
$$\lim_{n \to \infty} \frac{1}{v_{n-k}} \log A(n,k) = 0.$$

PROOF. First note that

$$\frac{1}{v_{n-k}} \log A(n,k)$$

$$= \frac{1}{v_{n-k}} \log \sum_{r=0}^{v_{n-k}-j_0^{n-k}-1} (f_k(e^{-\theta}))^r P(Z_{n-k} = r + j_0^{n-k})$$

$$\leq \frac{1}{v_{n-k}} \log \sum_{r \geq 0} (f_k(e^{-\theta}))^r,$$

and thus $\limsup_{n \to \infty} \frac{1}{v_{n-k}} \log A(n,k) = 0$. Furthermore,

$$\frac{1}{v_{n-k}} \log A(n,k)$$

$$= \frac{1}{v_{n-k}} \log \sum_{r=0}^{v_{n-k}-j_0^{n-k}-1} (f_k(e^{-\theta}))^r P(Z_{n-k} = r + j_0^{n-k})$$

$$\geq \frac{1}{v_{n-k}} \log P(Z_{n-k} = j_0^{n-k})$$

$$= \frac{n-k}{v_{n-k}} \log p_{j_0}.$$

Hence, $\liminf_{n\to\infty} \frac{1}{v_{n-k}} \log A(n,k) = 0$, concluding the proof of the claim. Also,

$$\frac{1}{v_{n-k}} \log(f_{n-k}(f_k(e^{-\theta}))(f_k(e^{-\theta}))^{-j_0^{n-k}})$$

$$= \frac{1}{v_{n-k}} \log f_{n-k}(f_k(e^{-\theta})) - \frac{j_0^{n-k}}{v_{n-k}} \log(f_k(e^{-\theta})).$$

Now using (see [13])

(136) $$\lim_{n\to\infty} \frac{1}{j_0^n} \log f_n(s) = G(s) \qquad \text{for } 0 < s < 1,$$

it follows that

$$\lim_{n\to\infty} \frac{1}{v_{n-k}} \log(f_{n-k}(f_k(e^{-\theta}))(f_k(e^{-\theta}))^{-j_0^{n-k}})$$

$$= bG(f_k(e^{-\theta})) - b \log f_k(e^{-\theta}). \qquad \square$$

Finally, using Claim 2, it follows that if $b > 0$, then

(137) $$\lim_{n\to\infty} \frac{1}{v_{n-k}} \log E(e^{-\theta N_n^1(k)}) = bG(f_k(e^{-\theta})).$$

The continuity of the limit follows from the continuity of $f(\cdot)$. If $b = 0$, then from (128) it follows that

(138) $$E(e^{\theta N_n^1(k)}) = (\delta_{n-k})^{-1} (f_k(e^{-\theta}))^{v_{n-k}} \sum_{r \geq 0} (f_k(e^{-\theta}))^r P(Z_{n-k} = r + v_{n-k}).$$

CLAIM 3.
$$\lim_{n\to\infty} \frac{1}{v_{n-k}} \log \sum_{r \geq 0} (f_k(e^{-\theta}))^r P(Z_{n-k} = r + v_{n-k}) = 0.$$

PROOF. It is easy to see that $\limsup_{n\to\infty} \frac{1}{v_{n-k}} \log \sum_{r \geq 0} (f_k(e^{-\theta}))^r P(Z_{n-k} = r + v_{n-k}) = 0$ using the trivial estimate $P(Z_{n-k} = r + v_{n-k}) \leq 1$. Also,

$$\frac{1}{v_{n-k}} \log \sum_{r \geq 0} (f_k(e^{-\theta}))^r P(Z_{n-k} = r + v_{n-k})$$

$$\geq \frac{1}{v_{n-k}} \log P(Z_{n-k} = v_{n-k})$$

$$= \frac{1}{v_{n-k}} \log(m^{n-k} P(Z_{n-k} = v_{n-k})) + \frac{n-k}{v_{n-k}} \log m.$$

Hence by Theorem 1, $\liminf_{n\to\infty} \frac{1}{v_{n-k}} \log \sum_{r\geq 0}(f_k(e^{-\theta}))^r P(Z_{n-k}=r+v_{n-k}) = 0$. Thus if $b=0$, then

$$\lim_{n\to\infty} \frac{1}{v_{n-k}} \log E(e^{\theta N_n^1(k)}) = \log f_k(e^{-\theta}). \tag{139}$$
□

Theorem 4 follows from Proposition 2. Finally, that the rate functions are good follows a similar line of proof as Theorem 3. □

REMARK 8. When $v_{n-k} = j_0^{n-k}$, then $b=1$ and $P(R_n \geq a | Z_{n-k} \geq j_0^{n-k})$ satisfies an LDP with rate function $-G(f_k(e^{-I(x)}))$, where $I(x)$ is the Legendre–Fenchel transform of $\log E(e^{\theta Z_1})$. However, since

$$P(R_n \geq a) = P(R_n \geq a | Z_{n-k} \geq j_0^{n-k}), \tag{140}$$

the distributions of $\{R_n : n \geq 1\}$ satisfy the LDP with the good rate function $-G(f_k(e^{-I(x)}))$. In comparison with Case $\alpha < \infty$, the rate function associated with the large deviations of $R_n$ is not a constant. This brings out yet another difference (in the probabilistic structure of the process) between Cases $\alpha < \infty$ and $\alpha = \infty$. It is also interesting to note that $b$ occurs additively in the Schröder case while it occurs multiplicatively in the Böttcher case.

**4. Concluding remarks.** In this paper we studied the local limit problem and the conditional large deviations of a general class of random variables indexed by branching processes. We brought to the fore the role played by the product $p_1 m$ (or $\gamma m$) in determining the asymptotic rates of $P(Z_n = v_n)$ and its impact on large deviations. The quantity $p_1 m$ can be viewed as a parameter that determines how "fast" the supercritical process is growing. If $p_1 m > 1$, then the process is growing slowly (since with $m$ thought of as fixed, $p_1$ is not very small) while if $p_1 m < 1$ the process is growing fast (since $p_1$ is very small). The "critical" case is $p_1 m = 1$.

One of the initial motivations for the topics in this paper was an interest in a version of the Gibbs conditioning principle for branching processes. In the simplest case of i.i.d. random variables the Gibbs principle determines the behavior of the individual members of a set of i.i.d. random variables, conditioned on the average of the whole ensemble.

In the branching context, the role of the average is played by $R_n$ and one is led to consider quantities such as

$$P((Z_{n_i}, \ldots, Z_{n_j}) \in \cdot | R_n \in A)$$

or

$$P((\xi_{n_i}, \ldots, \xi_{n_j}) \in \cdot | R_n \in A).$$

Since the random variable $R_n$ depends only on the present (two "present" generations) one could use the above formulation to make inferences on the history of the population based on the present. Preliminary calculations suggest interesting results in this direction, particularly in the multitype setting. In these calculations it became apparent that careful estimates on the behavior of $R_n$ in terms of the local limit estimates on $Z_n$ would be needed and this led to the results in the present paper.

Several questions arise from our work. The most interesting, from the large deviation perspective, seems to be when one replaces the indexing sequence by a "more general" sequence of random variables and allows dependences between the indexed and the indexing sequences. These kinds of problems occur naturally in sequential analysis, insurance and risk analysis areas. The authors are considering these generalizations and will report the results in future publications.

**Acknowledgments.** Part of this work was carried out while the second author was visiting the Center for Mathematical Sciences at the University of Wisconsin, Madison, during the years 1999 and 2000 and he thanks the Center for research facilities and travel support. The first author thanks the University of Georgia's Department of Statistics for support during this work.

## REFERENCES

[1] ATHREYA, K. B. (1994). Large deviation rates for branching processes I. Single type case. *Ann. Appl. Probab.* **4** 779–790. MR1284985
[2] ATHREYA, K. B. and NEY, P. E. (1970). The local limit theorem and some related aspects of super-critical branching processes. *Trans. Amer. Math. Soc.* **152** 233–251. MR268971
[3] ATHREYA, K. B. and NEY, P. E. (1972). *Branching Processes*. Springer, Berlin. MR373040
[4] ATHREYA, K. B. and VIDYASHANKAR, A. N. (1993). Large deviation results for branching processes. In *Stochastic Processes* (S. Cambanis, J. K. Ghosh, R. L. Karandikar and P. K. Sen, eds.) 7–12. Springer, New York. MR1427295
[5] BAHADUR, R. R. and RANGA RAO, R. (1960). On deviations of the sample mean. *Ann. Math. Statist.* **31** 1015–1027. MR117775
[6] BASAWA, I. V. (1981). Efficient conditional tests for mixture experiments with applications to the birth and branching processes. *Biometrika* **68** 153–164. MR614952
[7] BIGGINS, J. D. and BINGHAM, N. H. (1994). Large deviations in the supercritical branching process. *Adv. in Appl. Probab.* **25** 757–772. MR1241927
[8] BINGHAM, N. H. (1988). On the limit of a supercritical branching process. *J. Appl. Probab.* **25A** 215–228. MR974583
[9] BLAISDELL, B. E. (1985). A method of estimating from two-aligned present-day DNA sequences and their ancestral composition and subsequent rates of substitution, possibly different in the two lineages, corrected for multiple and parallel substitutions at the same site. *J. Mol. Evol.* **18** 225–239.

[10] Brown, W. M., Pager, E. M., Wang, A. and Wilson, A. C. (1982). Mitochondrial DNA sequences of primates: Tempo and mode of evolution. *J. Mol. Evol.* **18** 225–239.

[11] Dembo, A. and Zeitouni, O. (1998). *Large Deviations Techniques and Applications.* Springer, New York. MR1619036

[12] Dubuc, S. (1970). La fonction de Green d'un processus de Galton–Watson. *Studia Math.* **34** 69–87. MR260042

[13] Dubuc, S. (1971). La densitè de la loi-limite d'un processus en cascade expansif. *Z. Wahrsch. Verw. Gebiete* **19** 281–290. MR300353

[14] Dubuc, S. (1971). Problémes relatifs á l'itération de fonctions suggérés par les processus en cascade. *Ann. Inst. Fourier* **21** 171–251. MR297025

[15] Dubuc, S. and Senata, E. (1976). The local limit theorem for Galton–Watson process. *Ann. Probab.* **3** 490–496. MR405610

[16] Efron, B. and Hinkley, D. (1996). Assessing the accuracy of maximum likelihood estimator: Observed versus expected information. *Biometrika* **65** 457–487. MR521817

[17] Heyde, C. C. (1975). Remarks on efficiency in estimation for branching processes. *Biometrika* **62** 49–55. MR375695

[18] Heyde, C. C. (1977). An optimal property of maximum likelihood with application to branching process estimation. *Bull. Inst. Internat. Statist.* **47** 407–416. MR617588

[19] Joffe, A. and Waugh, W. A. O'N. (1982). Exact distributions of kin numbers in a Galton–Watson process. *J. Appl. Probab.* **19** 767–775. MR675140

[20] Joffe, A. and Waugh, W. A. O'N. (1985). Exact distributions of kin numbers in multitype Galton–Watson population. *J. Appl. Probab.* **22** 37–47. MR776886

[21] Karlin, S. and McGregor, J. (1968). Embeddability of discrete-time simple branching processes into continuous-time branching processes. *Trans. Amer. Math. Soc.* **132** 115–136. MR222966

[22] Karlin, S. and McGregor, J. (1968). Embedding iterates of analytic functions with two fixed points into continuous groups. *Trans. Amer. Math. Soc.* **132** 137–145. MR224790

[23] Kelly, C. (1994). A test for DNA evolutionary models. *Biometrics* **50** 653–664. MR1309311

[24] Ney, P. E. and Vidyashankar, A. N. (2003). Harmonic moments and large deviation rates for supercritical branching processes. *Ann. Appl. Probab.* **13** 475–489. MR1970272

[25] Pakes, A. G. (1975). Nonparametric estimation in the Galton–Watson process. *Math. Biosci.* **26** 1–18. MR405755

[26] Severini, T. A. (1996). Information and conditional inference. *J. Amer. Statist. Assoc.* **90** 1341–1346. MR1379476

[27] Severini, T. A. (2000). *Likelihood Methods in Statistics.* Oxford Univ. Press. MR1854870

[28] Sweeting, T. (1978). On efficient tests for branching processes. *Biometrika* **65** 123–127. MR494750

[29] Sweeting, T. (1986). Asymptotic conditional inference for the offspring mean of a supercritical Galton–Watson process. *Ann. Statist.* **14** 925–933. MR856798

[30] Taïb, Z. (1992). *Branching Processes and Neutral Evolution.* Springer, Berlin. MR1176317

[31] Zygmund, A. (1988). *Trigonometric Series*, 2nd ed. Cambridge Univ. Press. MR933759

Department of Mathematics  
University of Wisconsin  
Madison, Wisconsin 53706  
USA

Department of Statistics  
University of Georgia  
Athens, Georgia 30602  
USA  
e-mail: anand@stat.uga.edu